\documentstyle[12pt]{article}

\newcounter{const}

\newcommand{\argm}{\mathop{\rm argmax}\limits}

\begin{document}

 \large


\begin{center}

{\bf Exponential Bounds for Random Sums. }\\

\vspace{3mm}

{\sc Migdashiev B.M.$^b$, \  Ostrovsky E.I.$^a$*} \footnote{Corresponding Author.}
\footnote{ Department of Mathematics 
and Computer Science, Ben - Gurion University, Beer - Sheva, 
84105, Ben - Gurion street, 4;  P.O. Box 61, Israel. \\
E - mail: galaostr@cs.bgu.ac.il}\\

\vspace{3mm}

-------------------------------------------------------------------------------\\
{\footnotesize  {\bf Abstract.}  We construct a non - improved exponential bounds for
 distribution of normed sums of i.,i.d. random variables with random numbers of summand.}\par
{\it Key words:} Random sum, exponential estimation, Orlicz spaces, martingales,
saddle - point method.\\
--------------------------------------------------------------------------------\\

\end{center}

\vspace{2mm}

 {\bf 1. Introduction.}\par
 Let $ (\Omega,F,\mu) $ be a probability space, $ \{\xi(i) \}, i = 1,2,\ldots  $ be a 
sequence of independent  identical distributed (i; i.d) centered:  $ \ {\bf E} \xi(i) = 0 $ 
random variables (r.v) with finite non - trivial variance $ \sigma^2 = {\bf D}\xi(i) \in 
(0,\infty), $ and let $ \eta, \ \eta \ge 1 $ be an integer r.v. with finite first moment 
$ {\bf E} \eta = A, $ where $ A \in [2,\infty). $ We assume at first that the  r.v. 
$ \eta $  and the sequence $ \{\xi(i)\} $ are independent.\par
 We will denote for arbitrary r.v. $ \tau $ and $ x = const \ge 0 $ the tail function 
$$
T(\tau,x) = \max( {\bf P}(\tau \ge x), {\bf P}(\tau \le -x) ),
$$
 will write for the r.v. $ \xi(1) \ \  R(x) = T(\xi(1),x), $
and we define the so - called normed  random sum and corresponding uniform tail function 
$$
S = \sum_{i=1}^{\eta} \xi(i) /(\sigma \sqrt{A}), 
$$

$$
V(x) = V( Law(\eta), Law(\xi(i)), x) = \sup_{A \ge 2} T(S,x).
$$

 In the case if the r.v. $ \eta -1 $ has a Poisson distribution 
$ Pois(A) $ with parameter $ A: \ {\bf P}(\eta-1 = n) = A^n \exp(-A)/n!, \ n = 0,1,2,\ldots $ and 
$ Law(\eta-1)\in \cup_{A \ge 2} \{ Pois(A)\} $  we will write $ V(Pois; Law(\xi(i)), x) =
 V(\cup_{A \ge 2} \{Pois(A)\}, Law(\xi(i)), x); $
  if the r.v. $ \ \eta $ has a geometrical distribution 
$ G(A): \ {\bf P}(\eta = n) = A^{-1} (1-1/A)^{n-1}, n = 1,2,\ldots, \ A \ge 2$ and 
$ Law(\eta) \in \cup_{A \ge 2} \{G(A)\} $ we will write $ V(G; Law(\xi(i)), \ x)
= V(\cup_{A\ge 2} \{ G(A)\}, Law(\xi(i)),x); $ for the case 
of {\it all} distribution r.v. $ \eta $ under the  condition that
 $ \ {\bf E} \eta = A, \  \exists A \ge 2 $ we will use the notations correspondently
$ Dis(A), \ Dis = \cup_{A \ge 2} \{Dis(A)\} $ and 
$ V(Dis, Law(\xi(i)), x) = V(\cup_{A \ge 2} ( \{ Dis(A) \}, Law(\xi(i)),x).  $ \par 
 
{\bf Our goal is the bide - side exponential 
 estimating $ V(x) $ at $ x \to \infty, \  x \ge 2 $ in the terms of distributions }
 $ Law(\xi(i)), \ Law(\eta). $ \par 
 
 We have for all distribution $ \xi(i) $ under the conditions  $ {\bf E} \xi(i) = 0,\ 
{\bf D} \xi(i) \in (0, \infty),$ since $ {\bf E} \xi(i) = 0, \ {\bf D} \xi(i) \in (0,\infty): $
$$ 
V(Dis, Law(\xi(i)), x) \le \min \left(1, x^{-2} \right)) 
$$
by virtue of Chebyshev inequality; this estimation is called trivial.\par
 There are many publications 
about the moment estimations $ {\bf E }|S|^p $ and  statistical applications 
of those estimations (see, for example, (Gut A., 1988), (Gut A., 2003):
$$
|S|_p \le B(p) \ |\eta|_p^{1/2} \ |\xi(1)|_p, \ \ p \ge 2, \eqno(0)
$$
where $ B(p) $ is a constant in the famous Burkholder inequality for martingales.
Here and further for arbitrary r.v. $ \zeta $
$$
|\zeta|_p = {\bf E}^{1/p} |\zeta|^p = ||\zeta||L_p(\Omega, {\bf P}).
$$
It is proved in (Hitczenko, 1990) that the best boundary for $ B(p) $ at $ p \ge 2, \ 
p \to \infty $ is $ B(p) \le C \ p \ /\log p, $
$ C $ is an {\it absolute } constant.  Note that the estimation (0) is proved in the 
case if r.v.  $ \eta $ is the {\it stopping time } for the sequence $ \{\xi(i) \}, $
for example, if $ \forall n = 1,2,\ldots $ r.v.'s $ \xi(1), \xi(2), \ldots, \xi(n) $
and the event $ \{\eta = n \} $ are independent. \par 
 See also (Gine at al., 2003) etc. For non - random sums, i.e. in the case 
$ {\bf P} (\eta = n) = 1 $ 
for some $ n = 1,2,\ldots $ the exponential bounds for $ V(x) $ are constructed in 
(Buldygin at al., 1992), (Ostrovsky, 1999,
 p. 28); for more generally case if $ \{\xi(i) \} $ are the martingale -
differences the exponential bound for $ T(S,x) $ is derived in (Lesign at al., 2001). \par
 For reader convenience 
we recall here some results of (Buldygin at al, 1992),
(Lesign at al., 2001).
 Let us define for some tail function $ T(x), $ i.e. for the function $ T = T(x) $ under the
conditions:  $ T(0) = 1, \ T(\infty) = 0,$
monotonically non - increasing  and right continuous  with finite second moment 
$ |\int_0^{\infty} x^2 \ dT(x)| < \infty $ the operator 
$$
W[T](x) = \min \left(1, 4 \inf_{z > 0} \left[\exp(-x^2/(8z^2) - \int_z^{\infty} x^2 \ 
d T(x) \right] \right).
$$
 Further, put
$$
\varphi(\lambda) = \max_{\pm} \log {\bf E} \exp(\pm \lambda \xi(i)).
$$
 This definition is non - trivial only if the variable $ \xi(1) $ satisfies the so - called 
Kramer condition:
$$
\exists \lambda_0 \in (0,\infty), \ \forall \lambda \in (-\lambda_0, \lambda_0) \ \Rightarrow 
\varphi(\lambda) < \infty;
$$
in other case we set $ \varphi(\lambda) = \infty \ \forall \lambda \ne 0.$ \par 
 Let us introduce the function $ \chi(\lambda) = \sup_{n = 1,2,\ldots} 
n \varphi(\lambda/\sqrt{n}), $
$$
\chi^*(x) = \sup_{\lambda} (\lambda x - \chi(\lambda)),
$$
$$
Q[R](x) = \min \{W[R](x), \exp(-\chi^*(x)) \}.
$$
{\bf Lemma 1.} (Buldygin at al., 1992), (Lesign at al., 2001).
$$
\sup_n T \left(n^{-1/2} \sum_{i=1}^n \xi(i), \ x \right) \le Q[R ](x), 
\ \ x \ge 0. \eqno(1) 
$$
 Let us introduce the Orlicz spaces $ G(m,r) $ (in order to describe the examples) 
 of random variables as the set of all r.v. $ \{ \tau \} $  with finite norm
$$
||\tau||_{m,r} = \sup_{p \ge 2} |\tau|_p \ p^{-1/m} \ \log^{r/m} p.
$$
 Here $ m = const >0, \ r = const \in R^1. $ It is easy 
to verify that $ G(m,r) $ is isomorphic to the Orlicz space with $ N - $ function 
$$
N(u) = \exp \left(|u|^m \log^{r}|u| \right), \ \ |u| \ge 2,
$$
and that $ \tau \in G(m,r) $ if and only if 
$$
T(\tau,x) \le \exp \left( - C_1(m,r) x^m \log^r(C_2(m,r) + x) \right). \eqno(2)
$$
 See (Buldygin at al, 1992, p.351). \par 
For example, assume that  the r.v. $ \tau $ has Poisson distribution with 
parameter $ A; \ A \ge 2. $ 
Then for some non - trivial positive absolute constants $ C_1, C_2 $ and all $ p, x \ge 2 $
$$
|\tau - A|_p \le C_1 \ \sqrt{A} \ p/\log p,
$$
or 
$$
{\bf P}(|\tau - A|/\sqrt{A} > x) \le  \exp \left(-C_2 x \ \log x \right).
$$
 Let us suppose, for example, that for some $ m = const > 0, r \in R^1 \  \xi(1) \in G(m,r). $ 
We define the following  functions $ M = M(m,r), \ L = L(m,r): $ at $ m \in (0,1)  $
or $ m = 1, r < 0 \ \Rightarrow  \  M = 2m/(m+2), L = 2r/(m+2);  $ at 
$  m=1, r \ge 0 $  or $ m \in (1,2), r < 0  \ \Rightarrow  M = m,  L = r; $ 
 at $ m=2, r \ge 0 $  or $  m > 2 \ \Rightarrow M = 2, L = 0. $ We can define formally in 
the case $ m = +\infty, \  r \in R  \ \Rightarrow \  M = 2, \ L = 0. $ \par
 It follows from (1)
$$
  \ \sup_n T \left(\sum_{i=1}^n \xi(i)/\sqrt{n}, \ x \right) \le 
\exp \left(- C_3 \ x^M \ \log^L(C_2 + x)  \right), \eqno(3)
$$
$ C_{2,3} = C_{2,3}(m,r), \ $  or, equally, 
$$
\sup_{n \ge 1} ||\sum_{i=1}^n \xi(i)/\sqrt{n}||_{M,L} \le C(m,r) ||\xi(1)||_{m,r}. 
$$
 It is proved in (Ostrovsky, 1999, p. 34)  that in the case $ m > 1 $ the estimation (2) is exact 
at $ x \to \infty.$ \par
 {\it In this  paper, the letter $ C, C_j(\cdot) $ will denote positive finite
various non - essentially constant which may differ from one formula to the next and which 
does not  depend upon $ x,n.$ We make no attempt to obtain the best values for these constants.}\par

{\bf 2. Main result. Upper bound. Examples.} \par
{\bf Theorem 1.}
$$
V( Law(\eta), Law(\xi(i)), \eta,  x) \le {\bf E} \ Q[R]( \sigma \ x \sqrt{A/\eta}). \eqno(4)
$$
{\bf Proof.} We will assume without loss of generality $ \sigma = 1.$
 We receive from (1), using the formula of full probability and denoting 
$ q_n = q_n(A) = {\bf P}(\eta=n): \ \ V(x) = $ 
$$
 {\bf P} \left(\sum_{i=1}^{\eta} \xi(i) > x \sqrt{A} \right) = \sum_{n=1}^{\infty} q_n 
{\bf P} \left(\sum_{i=1}^n \xi(i)/\sqrt{n}  > x \sqrt{A/n} \right) \le 
$$
$$
\sum_n q_n \ Q[R](x\sqrt{A/n}) = {\bf E} \ Q[R](x \sqrt{A/\eta} ).
$$
{\it Example 1.} Let us suppose here that the r.v. $ \eta $  has a geometric distribution $ G(A) $
with parameter $ A, \ A \ge 2, $ and that $ \exists m > 0, \exists r \in R 
 \ \ L( \xi(1)) \in G(m,r). $ We assert that at $ x \ge 2 \ \Rightarrow
 \sup_{A \ge 2} V(G(A), G(m,r), x) \le $
$$
 C_3(m,r) \exp \left( - C(m,r) x^{2M/(M+2)} \ (\log x)^{2L/(M+2)} \right)
$$
$$
 \stackrel{def}{=}  J(M,L; x) = J(x). 
$$ 
{\bf Proof.}  We have using  theorem 1: 
$$
 V(G(A), G(m,r), x)  \le  \sum_{n=1}^{\infty} (A-1)^{-1} \exp ( n[\log(1-1/A)] - 
$$

$$
 C x^M \ n^{-M/2} A^{M/2} n^{-M/2} (\log^L ( C_2 +x \sqrt{A/n}) ) \le (A-1)^{-1}  \ \times
$$

$$
 C_1 \sum_{n=2}^{\infty} \exp \left(-C(n/A + x^M (A/n)^{M/2} 
(\log^L( C_2 + x \sqrt{A/n}) \right) \stackrel{def}{=}
$$

$$
  \sum_{n=2}^{\infty} a(n;A,x).
$$
 We denote

$$
N_0 = N_0(A,x) = \argm_{n \ge 2} \ a(n;A,x) = A \ n_0(x), 
$$
 where at $ x \to \infty \ \Rightarrow $

$$ 
N_0/A \sim C_2(M,L)  \ x^{2M/(M+2)} \left( \log x  \right)^{2L/(M+2) }.
$$

 Since the function $ n \to a(n;A,x) $ is monotonic in the  intervales $ [1, N_0] $ and
$ [N_0, \infty), $ we can estimate 

$$
 C_2 \ (A-1) \ V(G, G(m,r),x) \le 
$$
$$
 \int_1^{\infty} \exp \left(-C \left((y/A) + x^M (\sqrt{A/y})^M \
\log^L(C_2 + x \ \sqrt{A/y}) \right) \right) \ dy =
$$

$$
 C_3 A \  \int_{1/A}^{\infty} \exp \left( -C \left(z - x^M z^{-M/2}  \ \log^L(C_2 + x/\sqrt{z}) \
\right) \right) dz \le 
$$

$$
C_3 \ A \ \int_0^{\infty} \exp \left( - C(z - x^M z^{-M/2} \ \log^L(C_2 + x /\sqrt{z})
 \right) dz.
$$

 Let us denote $ \beta = 2M/(2+M), $

$$
R(x,v) = v + (2/M) \ v^{-M/2} \ \log^L \left(C_2 + x^{2/(M+2)} \ v^{-1/2} \right), 
$$

$$
U(x,v) = v + (2/M) v^{-M/2} \ \log^L x, \ x \ge 2,v > 0;
$$

$$
S(x,v) =  x^{\beta} \left(v + (2/M)  v^{-M/2} \ \log^L \left(C_2 + x^{2/(M+2)} v^{-1/2} \right) \right). 
$$

Let $ \Delta = \Delta(M,L) $ be the arbitrary function on $ M,L $ so that 
$ \Delta > 2|L|/(M+2). $
 After the substitution $ z = x^{\beta} \ v $ we receive:  $ \sup_{A \ge 2} 
V(G(A), G(m,r), C_1 x)/C_2 \le $

$$
C x^{\beta}\int_0^{\infty}\exp \left( - S(x,v) \right) dv \stackrel{def}{=} I(M,L;x) = I(x).
$$
 We have $ I = x^{\beta}(I_1 + I_2 + I_3), $ where
$$
I_1 = \int_0^{\log^{-\Delta} x}, \ I_2 = \int_{\log^{-\Delta} x}^{\log^{\Delta} x}, \ 
I_3 = \int_{\log^{\Delta} x}^{\infty} \exp(-S(x,v)) \ dv.
$$

Since at $ v \in (0, \log^{-\Delta} x ) \cup (\log^{-\Delta}x, \ \log^{+\Delta} x) )
 \ \Rightarrow $
$$
S(x,v) = x^{\beta} R(x,v) \ge C_3 U(C_4 x, v) \ge C x^{\beta} \log^{L + \Delta M/2} x,
$$
we conclude 

$$
 I_1 \le \log^{-\Delta} x \cdot \exp \left(-C x^{\beta} \ \log^{L + \Delta M/2} x 
\right) \le C_5 \  J(M,L; C_6 x); 
$$

$$
 x^{\beta} I_2 \le \left( \log^{\Delta} x  \right) \ \exp \left(  -C \min_{ \{v \in \log^{-\Delta}x, 
\log^{\Delta} x \} } U(x,v) \right) \le 
$$

$ C_7 J(M,N; C_8 x). $  Further,
$$
x^{\beta} \ I_3 \le \int_{\log^{\Delta} x}^{\infty} \exp \left( - x^{\beta} \ v \right) \ dv
 \le C_9 \ J(M,N; C_{10} x).
$$
 Following,  $ V(G,G(m,r),x) \le C_{11}(m,r) \times $
$$
 \exp \left( - C_{12}(m,r) x^{2M/(M+2)} \ \left(\log x \right)^{2L/(M+2)} 
 \right), \ x \ge 2.\eqno(5)
$$ 

 On the other hand 

$$
  \sup_{A \ge 2} \ || S||_{2M/(M+2), \ 2L/(M+2)} \le C_6 || \xi(1)||_{m,r}.
$$
{\it Example 2.} Let us now suppose again that  $ L\{\xi(i) \} \in G(m,r) $ for some $ m > 0,
r \in R $ and assume that the r.v.
$ \eta -1 $ has a {\it Poisson } distribution with parameter $ B = A -1; \ A \ge 2. $ 
 It follows  from theorem 1 and Stirling's formula that  $ V(Pois(A), G(m,r),x) \le $ 
$$
 C \sup_{B \ge 1} \sum_{n=1}^{\infty}  \exp ( -B + n \log B - \log n! -
$$

$$
  C x^M B^{M/2} n^{-M/2} (\log^L ( C_2 + x \sqrt{B/n}) ) ) \le
$$

$$
 C \sup_{B \ge 1} \sum_{n=1}^{\infty} \exp(-B  - n \log n +
$$

$$
 n - Cx^M B^{M/2} n^{-M/2} (\log^L(C_2 + x \sqrt{B/n} )) \le
$$
$$
 C \sum_{n=1}^{\infty} \sup_{B \ge 1} \exp(-B - n \log n + n -
$$

$$
 C x^M B^{M/2} n^{-M/2} (\log^L( C_2 + x \sqrt{B/n} )). \eqno(6)
$$

 It is easy to verify that the maximum of arbitrary member of the right - side  (6) 
over $ B \ge 1 $  for sufficiently greater values $ x \ge x_0, \ 
x_0 = const \ge 2 $  is achieved at $ B = 1. $ Therefore $ V(Pois(A), G(m,r), x)/C \le $

$$
\le \sum_{n=2}^{\infty} \exp \left(
 -1 - n \log n +n - C x^M n^{-M/2} \ \log^L(C_2 + x/\sqrt{n} ) \right).
$$
 Since for $ x \ge 2 $

$$
\sum_{n \ge  x^2} \exp(-n \log n + n) \le \exp \left( - C x^2 \log x \right),
$$
we have:

$$
V(Pois(A), \  G(m,r),\ x)/C \le \exp( - C x^2 \log x) + 
$$

$$
 \sum_{n \in [1, x^2]} \exp \left( n - n \log n - C x^M n^{-M/2} \ \log^L( C_2 + x /\sqrt{n}) \right). 
$$
 Let us denote $ N_1 = N_1(x) =  $
$$
 \argm_{n \in [1, x^2]} \left(n - n \log n - C x^M n^{-M/2} \ \log^L(C_2 + 
x/\sqrt{n} ) \right);
$$ 
it is easy to calculate:
$$
N_1(x) \asymp C(M,L) x^{2M/(M+2)} \ (\log x)^{(2L-2)/(M+2)}, \ \ x \to \infty.
$$
 We obtain  for values $ x\ge 2: $ 
$$ 
V(Pois, G(m,r), x)/C(m,r) \le \exp \left(-C x^2 \ \log x \right) +
$$
$$
 (x^2 +1)  \exp \left(N_1 - N_1 \log N_1 - C x^M N_1^{-M/2} \ \log^L(C_2+x/\sqrt{N_1}) \right)  \le 
$$
$$
 C_3 \exp \left( - C(m,r) x^{2M/(M+2)} \ (\log x)^{(2L + M)/(M+2)} \right). \eqno(7)
$$
 In the case $ m = \infty, $ i.e. if the variable $ \xi(1) $ is bounded $ ( mod \ \ {\bf P}), $
then $ 2M/(M+2) = 1, \ (2L+M)/(M+2) = 1/2. $ \par
  We can rewrite (7)  in the considered case $ Law(\eta -1) = Pois, Law(\xi(i) \in G(m,r) $
in the terms of $ G(m,r) $ spaces:

$$
\sup_{A \ge 2} ||S||_{2M/(M+2), \ (M+2L)/(M+2)} \le C(m,r) ||\xi(1)||_{m,r}.
$$

 In the case if $ Law\{ \xi(i) \} \in \cup_{m > 1} \{G(m,r)\} $ the estimation (7) improves  
some  result of (Gine at al, 2003). For instance, if $ m \in (1,\infty), \ r = 0, $ 
then from (7) it follows the inequality: $ \ p \ge 2 \ \Rightarrow $
$$
\sup_{A \ge 2} |S|_p \le C_1 \ p^{1/2 + 1/\min(m,2)}/\sqrt{\log p}, \ p \ge 2, 
$$
but we receive from (Gine at all, 2003):

$$
\sup_{A \ge 2} |S|_p \le C_2 \ p^{1+1/m} /\log p, 
$$
and we obtain from (Gut, 1988), (Gut, 2003):

$$
\sup_{A \ge 2} |S|_p \le C_3 \ p^{1+1/m} \ \log^{-1/2} p.
$$

\vspace{4mm}
{\bf 3. Low bounds.} We will prove further that our estimations are non - improved 
in general case, for example, even for normal distribution of values $ \{\xi(i)\} $. \par
{\bf Theorem 2.} {\it We assert that for all  values $ m > 1 $ and sufficiently larges }
 $ x, x \ge x_0 = const \ge 2: \ \ \  V(G, G(m,r), x)  \ge $ 
$$
C_6(m,r)
 \exp \left(-C_7(m,r) x^{2M/(M+2)} \ \left(\log x \right)^{2L/(M+2)} \right), \eqno(8)
$$

$$
V(Pois, G(m,r), x) \ge 
$$
$$
C_8(m,r) \exp \left(-C_9(m,r) x^{2M/(M+2)} \ (\log x)^{(2L+M)/(M+2)}  \right).\eqno(9)
$$
{\bf Proof }  is very simple. It is enought to prove the inequality (8); the proposition (9) 
is proved analogously. Let $ \xi(i) $ be independent symmetrically distributed r.v. with 
distributions

$$
{\bf P}(|\xi(i)| > x) = \exp \left(-x^m \ [\log(C(m,r) + x)]^r  \right), \ x \ge 0, 
$$
 and  let us introduce the even smooth  convex 
 function $ \varphi_{m,r}(\lambda) = \log \ {\bf E} \ \exp(\lambda \xi(i)), \ 
\lambda \in (-\infty,\infty). $ It is proved in (Buldygin at al., 1992, p.341), (Ostrovsky, 
1999, p.34) that
$$
C_1 \varphi_{m,r}(\lambda) \le \psi_{m,r}(\lambda) \le C_2 \varphi_{m,r}(\lambda),
$$
where $ \psi_{m,r}(\lambda) = \lambda^2, \lambda \in [-2, 2];$
$$
 \psi_{m,r}(\lambda) = C_3|\lambda|^{m/(m-1)} \ [\log(C_4+|\lambda|)]^{-r/(m-1)} 
$$
in the case $ |\lambda| > 2. $ Since the r.v. $ \{\xi(i) \} $ are i., i.d., we have for 
the non - random sum:
$$
\log {\bf E} \exp \left(\lambda n^{-1/2} \sum_{i=1}^n \xi(i) \right)= 
n \varphi_{m,r}(\lambda/\sqrt{n}) \asymp
$$

$  n \psi_{m,r}(\lambda/\sqrt{n}), $  where the symbol $ \asymp $ is understood 
uniformly on $ n; \ \lambda \in R: $ 

$$
0 < C_1(m,r) \le \inf_{n \ge 1} \ \inf_{\lambda \in R} \ \frac{n \varphi_{m,r}(\lambda/\sqrt{n}) }
{n \psi_{m,r}(\lambda/\sqrt{n}) } \le 
$$
$$
 \sup_{n \ge 1} \ \sup_{\lambda \in R} \ \frac{n \varphi_{m,r}(\lambda/\sqrt{n})}
{n \psi_{m,r}(\lambda/\sqrt{n})} \le C_2(m,r) < \infty. 
$$

 We conclude at $ x \ge 2, \ A \ge 2:$
$$
{\bf P} (S > x) = \sum_{n=1}^{\infty} A^{-1} (1-1/A)^{n-1} {\bf P} 
\left(\sum_{i=1}^n \xi(i)/\sqrt{n} >
x \sqrt{A/n} \right).
$$
 We deduce, choosing again in this sum only the member with  $ n = N_0 $  
(recall that $  N_0 = N_0(x)) $ and using the main result of paper (Bagdasarov at al, 1995):
$$
{\bf P} \left(n^{-1/2} \sum_{i=1}^n \xi(i) > u \right) \ge 
\exp \left(-C_1 u^M \ \log^L(C_2 + u) \right),
$$
where $ C_1, C_2 = C_1(m,r), C_2(m,r); \ \  u = u(n) \ge 2: $ 
$$
 {\bf P}(S > x) \ge (A-1)^{-1} (1-1/A)^{N_0} {\bf P} \left(\sum_{i=1}^{N_0} \xi(i) /\sqrt{N_0}
> x \sqrt{A/N_0}  \right) \ge
$$

$$
 C_6(m) \exp \left(-C_7(m) x^{2M/(M+2)} \ [\log x]^{2L/(M+2)} \right).  
$$

\vspace{3 mm}
{\bf Theorem 3.}  For all values $ x \ge 3 $ 
$$
V(Dis, N(0,1), x) \ge C \  x^{-2}. \eqno(10)
$$
 Here $ N(0,1) $ denotes the normal distribution with parameters $ 0, 1; $ 
and $ C $ is an absolute constant. \par
{\bf Proof.} We suppose now $ L(\xi(i)) = N(0,1), $ i.e. 
$$ 
{\bf P}(\xi(i) > x) = \Psi(x) = (2 \pi)^{-1/2} \int_x^{\infty} \exp(-y^2/2) dy. 
$$
 We define $ \alpha = 1/ Ent[x^2], \ Ent[z] $ denotes the integer part of $ z, $ for 
$ x \ge 3 $ and choose the r.v. $ \eta $ by the following way:
$ {\bf P}(\eta = 2) = 1- 1/\alpha, \ {\bf P}(\eta = 1/\alpha ) = \alpha.$
 Then $ A = {\bf E} \eta = 3 - 2 \alpha \ge 25/9 > 2;  $ 
$$
{\bf P}(S>x) > \alpha {\bf P} \left(\sum_{i=1}^{1/\alpha} \xi(i) > x \sqrt{A} \right) =
$$

$$
 (1/[x^2]) \  \Psi \left(x \sqrt{3-2/[x^2]} \sqrt{1/[x^2] } \right) \ge x^{-2} \ 
\Psi \left(3 \sqrt{3/8} \right) =   Cx^{-2}.
$$

\vspace{2mm}

{\bf 4. Upper exponential bound for stopping time. } In this section we  will obtain the 
exponential bounds for the tails of distribution r.v. $ S $ in the case if $ \eta $ is the 
{\it stopping  time } for the sequence $ \{ \xi(i) \}, $ in addition to the moment 
estimations of $ S $ in (Gut, 1988), (Gut, 2003).  Recall that again 
$  {\bf E} S = 0 $ and $ {\bf D} S = 1.$ \par
 
{\bf Theorem 4.} {\it Assume that the r.v. $ \eta $ belong to the space $ G(m,r) $ for 
some $ m > 0, \ r \in R:  \ \eta \in G(m,r)  $ and is the stopping time 
for the sequence $ \{ \xi(i) \}, $  where $ \xi(1) \in G(a,b), \ a = const >0, 
 b = const:  \ \forall x \ge 2  $

$$
T(\eta,x) \le \exp \left(-C \ x^m \ \log^r x \right), \ T(\eta,x) \le \exp \left(-C \ x^a 
\ \log^b x \right).
$$
Denote 
$$
q = \frac{2am}{2am + 2a +m}, \ w = \frac{2arm + mb + 2am}{2am + 2a +m}.
$$
 We assert that at } $ x \ge 2 $
$$
\sup_{A \ge 2} \sup_{\eta: {\bf E}\eta = A} T(S,x)  \le \exp \left( - C(a,b,m,r) x^q 
\ \log^w x  \right). \eqno(11)
$$
{\bf Proof of theorem 4. } It follows from our conditions and the theory of $ G \ - $ 
spaces (2) that for all $  \ p \ge 2 $

$$
|\xi(1)|_p \le C_1 p^{1/m} \ \log^{-r/m} p, \  \ |\eta|_p \le C_2 p^{1/a} \ \log^{-b/a} p.
$$

 We obtain using the inequality (0) with optimal constant $ B(p): $

$$
|S|_p \le C_3 \ p^{1 + (1/2a) + 1/m} \ \log^{-1 - (b/(2a) - (r/m) } p = C_3 p^{1/q} 
\ \log^{ - w/q} p.
$$
 The last inequality is equivalent to (11), see (2).\\

\newpage

\vspace{3mm}
{\bf References} \\
\vspace{2mm}

Bagdasarov D.R., Ostrovsky E.I., 1995. Reversion of Chebyshev's Inequality. Probab. Theory 
Appl., v. 40 {\bf 4}, 737 - 742.\\

Buldygin V.V., Mushtary D.I., Ostrovsky E.I., Puchalsky M.I., 1992. New Trends in Probability 
Theory and Statistic. Mokslas; Amsterdam, New York.\\

Gine E., Mason D.M., Zaitsev A.Yu., 2003. The $ L_1 - $ norm Density Estimation Process. Annals 
Probab., v. 31, {\bf 2}, 719 - 768.\\

Gut A., 1988. Stopped Random Walks. Springer Verlag, Berlin - Heidelberg - New York.\\

Gut A., 2003. On the Moment Problem for Random Sums. Journal of Appl. Probab., v. 40 $ N^o $ 3, 
707 - 802.\\

Hitczenko P., 1990. Best constant in martingale version of Rozental's Inequality. 1990. Annals 
Probab., v. 18,  1656 - 1668.\\

Lesign Emm., Dalibor, 2001. Large deviations for martingales. Stochastic Process. 
Appl;., {\bf 96,} 143 - 159.\\
 
Ostrovsky E.I., 1999. Exponential Estimations for the Random Fields. OINPE, Moscow (in Russian).\\

 \newpage
Ben - Gurion University, Department of Mathematic.\\

\vspace{3mm}
ISRAEL, Beer - Sheva city,  84105, Ben - Gurion street, 4. \  P.O.Box 61.\\

\vspace{3mm}
e - mail: \  Galaostr@cs.bgu.ac.il\\

\end{document}